# A Simple Solution to a Major Problem: Proof of the Riemann Hypothesis

Fayang G. Qiu
Email: qiu_fayang@gibh.ac.cn
Laboratory of Molecular Engineering
188 Kaiyuan Blvd., 5th Floor, Building D
The Science Park of Guangzhou, 510530, China

**Abstract**: Based on Riemann's symmetrical reflection functional equation of the zeta function, we have found that the $\sigma$ values satisfying $\zeta(s) = 0$ must also satisfy both $|\zeta(s)| = |\zeta(1 - s)|$ and $|\Gamma(s/2)\zeta(s)| = |\Gamma((1 - s)/2)\zeta(1 - s)|$. We have shown that $\sigma = ½$ is the only numeric solution that satisfies this requirement.

## Introduction

The Riemann zeta function expressed by the following series [1]

$$\zeta(s) = \sum_{n=1}^{\infty} \frac{1}{n^s} = 1 + \frac{1}{2^s} + \frac{1}{3^s} + \ldots \qquad (1)$$

converges and defines a function analytic in the region {$s$ in **C**: Re($s$) > 1}. Since the series diverges in the other region {$s$ in **C**: Re($s$) ≤ 1}, (1) is no longer valid. Riemann realized that the zeta function can be redefined using an analytic continuation and he extended it in a unique way to a meromorphic zeta function $\zeta(s)$ defined for all complex numbers $s$ except $s = 1$, where the zeta function has a singularity. Riemann has also established[1] that the zeta function satisfies the functional equation:

$$\zeta(s) = 2^s \pi^{s-1} \sin\left(\frac{\pi s}{2}\right)\Gamma(1-s)\zeta(1-s) \qquad (2)$$

where $\Gamma(s)$ is the gamma function, which is an equality of meromorphic functions valid on the entire complex plane. Owing to the properties of the sine function, the functional equation implies that $\zeta(s)$ has simple zeroes at negative even integers, which are known as the trivial zeros. In addition to these trivial zeros, the extended zeta function also has many nontrivial zeros in the strip {$s \in$ C: 0 < Re($s$) < 1}, which is called the critical strip.

The nontrivial zeros are more important because not only their distribution is far less understood, but their studies have yielded important results concerning prime numbers and related subjects in number theory.[2] When Riemann calculated a few nontrivial zeros, it appeared that all had $\boldsymbol{R}[s] = ½$. Thus he hypothesized that all nontrivial zeros of the zeta function have real part $\boldsymbol{R}[s] = ½$, which is confirmed to be true for the first $10^{12}$ roots.[2] However, at present, the Riemann Hypothesis still remains to be proven. In this communiccation, we provide a straightforward yet simple solution to this major probelm.

## The Proof

Riemann has established the following the symmetric reflection functional equation[1]

$$\pi^{-\frac{s}{2}}\Gamma(\frac{s}{2})\zeta(s)=\pi^{-\frac{1-s}{2}}\Gamma(\frac{1-s}{2})\zeta(1-s) \tag{3}$$

where $s = \sigma + ib$. Taking the absolute value for both sides, one obtains

$$|\pi^{-\frac{s}{2}}\Gamma(\frac{s}{2})||\zeta(s)|=|\pi^{-\frac{1-s}{2}}\Gamma(\frac{1-s}{2})||\zeta(1-s)| \tag{4}$$

Since both $|\pi^{-s/2}\Gamma(s/2)|>0$ and $|\pi^{-(1-s)/2}\Gamma((1-s)/2)|>0$ are true for all $s$ in the complex plane, if $|\zeta(s)| \neq |\zeta(1 - s)|$, then, at least one of them is nonzero. According to (4), if $|\zeta(s)| > 0$, then $|\zeta(1 - s)| > 0$, and *vice versa*. In this case, neither $|\zeta(s)|$ nor $|\zeta(1 - s)|$ can be zero. Hence no nontrivial zeros of $\zeta(s)$ will occur in the area where $|\zeta(s)| \neq |\zeta(1 - s)|$. Consequently, all nontrivial zeros will be bound to the $\sigma$ values that satisfy:

$$\left|\zeta(s)\right|=\left|\zeta(1-s)\right| \tag{5}$$

It follows immediately that, if $|\Gamma(s/2)\,\zeta(s)| \neq |\Gamma((1-s)/2)\,\zeta(1 - s)|$, then, at least one of them is nonzero. According to (4), if $|\Gamma(s/2)\,\zeta(s)| > 0$, then $|\Gamma((1-s)/2)\,\zeta(1 - s)| > 0$, and *vice versa*. In this case, neither $|\Gamma(s/2)\,\zeta(s)|$ nor $|\Gamma((1-s)/2)\,\zeta(1 - s)|$ can be zero. Hence no nontrivial zeros of $\zeta(s)$ will occur in the area where $|\Gamma(s/2)\,\zeta(s)| \neq |\Gamma((1-s)/2)\,\zeta(1 - s)|$. Thus, the $\sigma$ values of all nontrivial zeros will satisfy:

$$\left|\Gamma\left(\frac{\sigma+ib}{2}\right)\zeta(s)\right|=\left|\Gamma\left(\frac{1-\sigma-ib}{2}\right)\zeta(1-s)\right| \tag{6}$$

After (4) and (6) are combined, one obtains (7):

$$\left(|\pi^{-\frac{\sigma+ib}{2}}|-|\pi^{-\frac{1-\sigma-ib}{2}}|\right)\left|\Gamma(\frac{\sigma+ib}{2})\zeta(\sigma+ib)\right|=0 \tag{7}$$

For all $|\zeta(\sigma + ib)| > 0$, including $|\zeta(\sigma + ib)| \to 0$, (7) holds if and only if (8) is true:

$$\left|\pi^{-\frac{\sigma+ib}{2}}\right|=\left|\pi^{-\frac{1-\sigma-ib}{2}}\right| \tag{8}$$

In this case, **(6) and (8) are mutually sufficient and necessary conditions to each other**. It is easy to see that $\boldsymbol{\sigma = ½}$ is the only numeric value that satisfies (8).

Consider the case (A - B) x = 0. Since x is continuous, it comes immediately that A = B for all $|x| \geq 0$. Similarly, Riemann's $\zeta(s)$ is continuous on the entire complex plane excepting $s = 1$, and the $\Gamma(s)$ function is continuous in the critical strip, the $\sigma$ values satisfying (7) for all $|\zeta(\sigma + ib)| \geq 0$ will lead to (8).[3] Therefore, $\boldsymbol{\sigma = ½}$ is a necessary condition for the Riemann zeta function to have nontrivial zeros, which is equivalent to the statement that **all nontrivial zeros of the zeta function have real part $\boldsymbol{R[s] = ½}$.**

In addition, since this proof is independent of the form of the zeta function, thus, **all the nontrivial zeros of the general *L*-functions will have real part $\boldsymbol{\sigma = ½}$.**